\documentclass[12pt]{article}
\usepackage{amsmath}
\usepackage{amsthm}
\usepackage{amssymb}
\usepackage{pdflscape}
\usepackage{geometry}
\usepackage{cleveref}
\usepackage[T1]{fontenc}
\usepackage[sc]{mathpazo}
\usepackage{enumerate}
\usepackage{etoolbox}
\usepackage{hhline}

\patchcmd{\section}{\scshape}{\bfseries}{}{}

\newcommand{\intersect}{\cap}

\newcommand{\isom}{\cong}

\renewcommand{\subset}{\subseteq}
\newcommand{\union}{\cup}
\newcommand{\Union}{\bigcup}

\newcommand{\F}{{\mathbb{F}}}
\newcommand{\Z}{{\mathbb{Z}}}
\newcommand{\K}{{\mathbb{K}}}
\renewcommand{\P}{{\mathbb{P}}}

\newcommand{\directSum}{\oplus}
\newcommand{\epimorphism}{\twoheadrightarrow}

\newcommand{\ijkl}{i+j + k + \ell}

\renewcommand{\epsilon}{\varepsilon}

\renewcommand{\phi}{\varphi}

\DeclareMathOperator{\Cay}{Cay}

\newcommand{\axiom}[1]{\mathsf{#1}}

\newcommand{\AxiomK}{\axiom{CCA}}

\numberwithin{table}{section}

\newtheorem{theorem}[table]{Theorem}
\newtheorem{proposition}[table]{Proposition}
\newtheorem{lemma}[table]{Lemma}
\newtheorem{corollary}[table]{Corollary}
\newtheorem{example}[table]{Example}

\theoremstyle{definition}
\newtheorem{definition}[table]{Definition}

\theoremstyle{remark}
\newtheorem{remark}[table]{Remark}

\newtheorem{question}[table]{Question}

\begin{document}
\title{Cayley graphs on elementary abelian groups of extreme degree have complete cores}
\author{Guang Rao\footnote{Email: \texttt{generalrao@hotmail.com} 
} \and Colin Tan}
\date{30 Jan 2025}

%\address
%{Claus Scheiderer,
%Fachbereich Mathematik und Statistik,
%Universt\"at Konstanz,
%Konstanz 78457,
%Germany}
%\email{claus.scheiderer@uni-konstanz.de}

%\address
%{Colin Tan,
%Department of Statistics \& Applied Probability,
%National University of Singapore,
%Block S16,
%6 Science Drive 2,
%Singapore 117546}
%\email{statwc@nus.edu.sg}

\maketitle
\numberwithin{equation}{section}

\begin{abstract}
Ne\v{s}et\v{r}il and \v{S}\'{a}mal asked whether every cubelike graph has a cubelike core. Man\v{c}inska, Pivotto, Roberson and Royle answered this question in the affirmative for cubelike graphs whose core has at most $32$ vertices. When the core of a cubelike graph has at most $16$ vertices, they gave a list of these cores, from which it follows that every cubelike graph with degree strictly less than $5$ has a complete core. We prove the following extension: if the degree of a cubelike graph is either strictly less than $5$ or at least $5$ less than the number of its vertices, then its core is complete and induced by a $\mathbb{F}_2$-vector subspace of its vertices. Thus we also answer Ne\v{s}et\v{r}il and \v{S}\'{a}mal's question in the affirmative for cubelike graphs with degree at least $5$ less than the number of vertices. Our result is sharp as the $5$-regular folded $5$-cube and its graph complement are both non-complete cubelike graph cores. We also prove analogous results for Cayley graphs on elementary abelian $p$-groups for odd primes $p$.
\end{abstract}

\section{Introduction} \label{sec: intro}

A Cayley graph $\Cay(A, C)$ on an abelian group $A$ (always written additively) is specified by its \emph{connection set} $C \subset A\setminus \{0\}$
        satisfying $-C = C$.       
Explicitly, two vertices $x, y \in V(\Cay(A, C)) := A$ are adjacent in $\Cay(A, C)$ % by definition 
if and only if $x - y \in C$.
Note that the conditions $-C = C$ and  $0 \notin C$ ensure that 
%the adjacency relation is symmetric (resp.\ irreflexive) so that
$\Cay(A, C)$ is undirected and loopless.
When $A = {\Z_2}^d$ is an elementary abelian $2$-group,
    such Cayley graphs are known variously as \lq\lq \emph{cubelike graphs}\rq\rq\ \cite{Lovasz1975} or \lq\lq \emph{binary Cayley graphs}\rq\rq\ \cite{BNT2015}.
More generally, we shall refer to Cayley graphs on elementary abelian $p$-groups, where $p$ is a prime, as \lq\lq\emph{$p$-ary Cayley graphs}\rq\rq. 

The following question was first asked for cubelike graphs by Ne\v{s}et\v{r}il and \v{S}\'{a}mal \cite{NesetrilSamal2008}
    and was raised by Ong for general $p$ in her undergraduate capstone project supervised by the first author \cite{Ong}:
\begin{question} \label{qn: conjecture}
For a given prime $p$,
    is the core of a $p$-ary Cayley graph necessarily also a $p$-ary Cayley graph?
\end{question}

%(All graphs in this paper are simple and undirected.) 
Recall that a graph is a \emph{core} if all its endomorphisms are invertible.
A \emph{core} of a graph $X$ 
    is a subgraph that is both a retract of $X$ and a core.
All cores of a given graph are isomorphic. 
The core of a graph $X$ is denoted by $X^\bullet$.

Man\v{c}inska, Pivotto, Roberson and Royle proved that, for every cubelike graph $X$, if its core $X^\bullet$ has at most $2^5 = 32$ vertices, then $X^\bullet$ is necessarily also cubelike \cite[Theorem 8.2]{MPRR2020}. Thus they partially answered Question \ref{qn: conjecture} in the affirmative for cubelike graphs with small cores.
If $X^\bullet$ has at most $2^4 = 16$ vertices, then they showed more specifically that $X^\bullet$ is either complete, or the halved $5$-cube $\frac{1}{2}Q_5 = \Cay({\Z_2}^4, \{i, j, k, \ell, i +j, i + k, i + \ell, j + k, j + \ell, k + \ell\})$, or the folded $5$-cube $\overline{\frac{1}{2}Q_5} \isom \Cay({\Z_2}^4, \{i, j, k, \ell, i + j + k + \ell\})$ \cite[Theorem 8.1]{MPRR2020}.
Here $\overline{X}$ denotes the complement of a graph $X$. 
We use the consecutive letters $i := e_1$, $j := e_2$, $k := e_3$, $\ell := e_4$, \dots, 
    to denote the first few vectors in an $\Z_p$-linear basis $\{e_1, \dots, e_d\}$ of ${\Z_p}^d$.
Below, we will often regard ${\Z_p}^d = ({\F_p}^d, +)$ as the underlying additive abelian group of the $d$-dimensional vector space ${\F_p}^d$ over the finite field $\F_p$ of order $p$.
    
This latter result implies the following:
\begin{theorem} \label{thm: consequenceOfMPRR}
Every cubelike graph $X$ with $\deg(X) < 5$ has a complete core.  
\end{theorem}

Here $\deg(X)$ denotes the common degree of all the vertices of a $p$-ary Cayley graph $X$. 
%Recall that Cayley graphs are vertex transitive. 

\begin{proof}
%Suppose that $X$ is a cubelike graph with $\deg(X) < 5$.
%We may assume without loss of generality that $X$ is connected, since the core of a cubelike graph is the core of any of its connected components.  
Assume without loss of generality that $X$ is connected.
Then $X$ has at most $2^{5 - 1} = 16$ vertices.
By %a theorem of Man\v{c}inska, Pivotto, Roberson and Royle 
\cite[Theorem 8.1]{MPRR2020}, $X^\bullet$ is either complete,  $\frac{1}{2}Q_5$ or $\overline{\frac{1}{2}Q_5}$.
But $\deg(\frac{1}{2}Q_5) = 10$ and $\deg(\overline{\frac{1}{2}Q_5}) = 5$, so $\deg(X^\bullet) \le \deg(X) < 5$ implies that $X^\bullet$ is complete.
\end{proof}
Our main result extends Theorem \ref{thm: consequenceOfMPRR} to $p$-ary Cayley graphs of either low or high degree.
%For clarity, we state our results for the cases of $p = 2$ and $p = 3$ separately.
Define the following function $\kappa$ on the set of all primes:
\begin{equation} \label{eq: counterexampleDegree}
\kappa(p) := \begin{cases}
5 & \text{if } p = 2; \\
12 & \text{if } p = 3; \\
2 & \text{if } p \ge 5 \text{ is prime}.
\end{cases}
\end{equation}
%Let $\elem{d} = (\Z/(2\Z))^{\directSum d}$ denote the additive elementary abelian $2$-group of order $2^d$, where $d$ is a nonnegative integer.
\begin{theorem} \label{thm: completeCores}
Let $p$ be a prime and let $d \ge 0$ be an integer. 
For every $p$-ary Cayley graph $X$ on $({\F_p}^d, +)$,
    if either $\deg(X) < \kappa(p)$ or $\deg(X) \ge p^{d} - \kappa(p)$, 
    then there is a $\F_p$-vector subspace $V \subset {\F_p}^d$ that induces a complete core of $X$. 
\end{theorem}

We prove Theorem \ref{thm: completeCores} in \S\ref{sec: proofOfTheorem}. This theorem leads immediately to the following partial affirmative answer to Question \ref{qn: conjecture} for $p$-ary Cayley graphs of extreme degree: 

\begin{corollary} \label{cor: conjectureForExtremeDegree}
Let $p$ be a prime.
Then every $p$-ary Cayley graph $X$ with
    either $\deg(X) < \kappa(p)$ or $\deg(X) \ge |V(X)| - \kappa(p)$
has a $p$-ary Cayley graph as its core.
\end{corollary}

\begin{proof}
By Theorem \ref{thm: completeCores}, $X^\bullet = X|_V = K_V= \Cay(V, V \setminus \{0\})$.
\end{proof}

Here $X|_V$ denotes the subgraph of a graph $X$ that is induced by a subcollection $V \subset V(X)$ of the vertex set of $X$. 

In fact, we obtain a full classification of the dimension of $V$ such that $X^\bullet = X|_V$
    given a $p$-ary Cayley graph $X$ 
        satisfying the degree bounds
            in Theorem \ref{thm: completeCores} (see Table \ref{tab: p=2ExtremeCardinality}, Tables \ref{tab: p=3LowCardinality}, \ref{tab: p=3HighCardinality}, and Table \ref{tab: p>=5ExtremeCardinality} below for the cases $p = 2$, $p = 3$, and $p \ge 5$ respectively).
In particular, 
    we show that if $\deg(X) \ge p^d - \kappa(p)$, then $\dim_{\F_p}(V) \ge d - 2$, so $|V(X^\bullet)| \ge p^{d - 2}$.
Thus the number of vertices of $X^\bullet$ grows arbitrarily large as $d \to \infty$,
    thus, even when $p = 2$, Theorem \ref{thm: completeCores} addresses cases not covered by \cite[Theorem 8.2]{MPRR2020}.    
%        its core is induced by some vector subspace $V \subset {\F_2}^d$ of dimension at least $d - 2$.
%For example, if $\deg(X) = 2^d - 1$ is maximum, then $X = K_{2^d}$ is the complete graph on $2^d$ vertices. Every endomorphism of $X$ is bijective and hence invertible, thus $X^\bullet = X = K_{2^d}$. 
%We give another example of Theorem \ref{thm: CCCForHighDegree} that will follow from our results below. 
%\begin{example} \label{eg: directionsetHasCardinality1}
%Let $d \ge 1$ 
%    and let $X_d$ be the cubelike graph with $2^d$ vertices having degree $2^d - 2$.
%Then $(X_d)^\bullet = K_{2^{d - 1}}$. 
%Explicitly, write $X_d = \Cay(({\F_2}^d, + ), C)$, where its connection set $C \subset {\F_2}^d \setminus \{0\}$ has cardinality $|C| = \deg(X_d) = 2^d - 2$. Then ${\F_2}^d \setminus C = \{0, v\}$ for some nonzero $v \in {\F_2}^d$ and is necessarily a $1$-dimensional $\F_2$-linear subspace of ${\F_2}^d$. Then in fact any linear hyperplane $H \subset {\F_2}^d$ which decomposes ${\F_2}^d = H \directSum ({\F_2}^d \setminus C)$ into a direct sum induces a complete core of $X_d$.   
%\end{example}
%Thus, for any infinite sequence $(X_d)_{d = 1}^\infty$ 
%    whose $d$-th term is a cubelike graph $X_d$ with $2^d$ vertices and $\deg(X_d) \ge 2^d - 5$,
%each of whose core ${X_d}^\bullet = K_{2^{d-1}}$ is a proper cubelike subgraph, where 
%the number of vertices of ${X_d}^\bullet$ is at least $2^{d - 2}$
%    and thus approaches infinity as $d \to \infty$. 
%So we obtain a partial affirmative answer to Question B by showing that there are cubelike graphs with arbitrary large cubelike cores. 

The degree bounds given in Theorem \ref{thm: completeCores} are sharp as sufficient conditions for $X^\bullet$ to be complete as witnessed by the following counterexample:
%    that ensure a $p$-ary Cayley graph has a complete core
%, as the following counterexample shows.    
%Recall the function $\kappa$ defined in \eqref{eq: counterexampleDegree}.
\begin{example} \label{eg: counterexamples}
For each prime $p$,
    there exists a $p$-ary Cayley graph $X_\ast$
            of degree $\kappa(p)$ 
    such that both $X_\ast$ and %its graph complement
    $\overline{X_\ast}$ are non-complete cores.
In particular, $X_\ast$ (resp.\ $\overline{X_\ast}$) is a $p$-ary Cayley graph of degree $\kappa(p)$ (resp.\ $|V(\overline{X_\ast})| - \kappa(p) - 1)$ whose core is not complete.
\end{example}
The details of Example \ref{eg: counterexamples} are given in \S\ref{sec: counterexamples}. 

%To the authors' knowledge, this paper is the first 
Other than \cite{MPRR2020},
    only Ong \cite{Ong} and Rotheram \cite{Rotheram2013} have written about Question \ref{qn: conjecture}. 
Ong computed the core of small ternary Cayley graphs
    and the core of $p$-ary Cayley graphs $X$ with $\deg(X) \le |V(X)| - 3$ for general $p$.
For $p\neq 3$, Rotheram showed that if the automorphism group of $X^\bullet$ is primitive in its action on $V(X^\bullet)$, then $X^\bullet$ is necessarily a $p$-ary Cayley graph. 

%\paragraph{Organisation of the rest of this paper.}
In \S\ref{sec: axiomatics}, we formulate an abstract property on a connection set that is sufficient for the core of its associated Cayley graph to be complete. We prove Theorem \ref{thm: completeCores} in \S\ref{sec: proofOfTheorem} and give details of Example \ref{eg: counterexamples} in \S\ref{sec: counterexamples}.    
We conclude by raising an open question. 

\section{A sufficient condition for a Cayley graph on an abelian group to have a complete core} \label{sec: axiomatics}

We give a sufficient condition on the connection set $C$ of a Cayley graph on an abelian group $A$
    for $\Cay(A, C)^\bullet$ to be complete.    
%All abelian groups are written additively.  
%Recall from \S\ref{sec: mainResult} that a subset $C \subset A \setminus \{0\}$
%    is said to be a \emph{connection set}
%        if $-C = C$.
Here and below, the symbol \lq$C$\rq\ shall be used with further mention to refer to a connection set.
%, which shall be understood to satisfy $-C = C$.       
Recall that subgroups $G, H \le A$ decompose $A$ into a \emph{direct sum} $A = G \directSum H$ if every $x \in A$ can be written uniquely as $x = g + h$ for some $g \in G, h \in H$. 

\begin{proposition} \label{prop: applicationOfAxiom}
Let $C \subset A \setminus \{0\}$, where $A$ is an abelian group.
If there are subgroups $G, H \le A$
    that decompose $A = G \directSum H$ such that
\begin{align}
G &\subset C \union \{0\}  \label{eq: AxiomKPartOne} \\
\text{and }H \intersect C &= \emptyset, \label{eq: AxiomKPartTwo}
\end{align}        
then the projection $p_G: \Cay(A, C) \to \Cay(A, C)|_G$
    from $A$ onto its direct summand $G$
    is a retraction onto a complete core of $\Cay(A, C)$. 
\end{proposition}

\begin{proof}
Let $X := \Cay(A, C)$ 
    and suppose that $C \subset A \setminus \{0\}$ 
        satisfies \eqref{eq: AxiomKPartOne} and \eqref{eq: AxiomKPartTwo} 
            for some direct-sum decomposition $A = G \directSum H$. 
To see that $p_G$ is a homomorphism,
    let adjacent vertices $x, x' \in V(X) = A$ be given. 
Then, by definition, $x - x' \in C$. 
Write $x = g + h$ and $x' = g' + h'$, 
%in terms of the direct-sum decomposition $A = G \directSum H$, 
where $g, g' \in G$ and $h, h' \in H$. Then
$
(g - g') + (h - h') = x - x' \in C
$.
Thus $g - g' \neq 0$, using \eqref{eq: AxiomKPartTwo}. 
Therefore $p_G(x) - p_G(x') = g - g' \in G \setminus \{0\} \subset C$, using \eqref{eq: AxiomKPartOne},
    so that $p_G(x)$ and $p_G(x')$ are adjacent in $X|_G$.

Since $p_G$ restricts to the identity map on $G$,
    thus $p_G$ is a retraction onto $X|_G$. 
Moreover, \eqref{eq: AxiomKPartOne} implies that the induced subgraph $X|_G = K_G$ is complete and hence a core. 
Therefore, $X|_G = X^\bullet$ is a complete core of $X$. 
\end{proof}

As this abstract property will be used very often,
    we make the following definition.
Let $A$ be an abelian group.

\begin{definition}
Say that $C \subset A \setminus \{0\}$
    satisfies the  
    $\mathsf{Complete\ Core\ Axiom}$ with respect to $(G, H)$ (or just $\AxiomK(G, H)$ or even $\AxiomK$ for short)
        if $G, H \le A$ are subgroups
            such that $A = G \directSum H$
                and both \eqref{eq: AxiomKPartOne} and \eqref{eq: AxiomKPartTwo} hold.

%Sometimes we omit mention of $(G, H)$ and say vaguely that $C$ satisfies $\AxiomK$. 
\end{definition}

Given a subset $S$ of a vector space $\K^d$ over a field $\K$,
   let $\langle S \rangle_{\K}$ (or just $\langle S \rangle$ if the base field is understood) denote the $\K$-linear span of $S$ in $\K^d$.
If $S$ is finite, say $S = \{v_1, \dots, v_t\}$,
    we also write $\langle S \rangle = \langle v_1, \dots, v_t\rangle$.

\begin{example}
Consider $X = \Cay(({\F_p}^d, +), C)$, where $C = \Union_{i = 1}^d (\langle e_i \rangle_{\F_p} \setminus \{0\})$. 
As an iterated Cartesian product, $X = \Cay((\F_p, +), \F_p \setminus \{0\})^d = {K_p}^d$.
For instance, when $p = 2$, then $X = \square^d$ is the $d$-dimensional hypercube.
Since $X$ satisfies $\AxiomK(\langle e_1 \rangle_{\F_p}, \{\lambda_1 e_1 + \cdots + \lambda_d e_d \in {\F_p}^d : \, \lambda_1 + \cdots + \lambda_d = 0\})$,
by Proposition \ref{prop: applicationOfAxiom}, $X|_{\langle e_1 \rangle} = K_p$ is a complete core of $X$.
\end{example}

%Here are some basic lemmas regarding $\AxiomK$.

\begin{lemma} \label{lem: complementAlsoSatisfiesAxiomD}
Let $C \subset A \setminus \{0\}$.
%, where $A$ is an abelian group.
If $C$ satisfies $\AxiomK(G, H)$, 
%for some directsum decomposition $A = G \directSum H$ of abelian groups,
then its complement $\overline{C} \subset A \setminus \{0\}$ 
    satisfies $\AxiomK(H, G)$.
\end{lemma}

\begin{proof}
Suppose $C \subset A \setminus \{0\}$ satisfies $\AxiomK(G, H)$.
%for some directsum decomposition $A = G \directSum H$ of abelian groups. 
Then of course $A = G \directSum H \isom H \directSum G$ also.
From \eqref{eq: AxiomKPartTwo}, 
$
H \subset A \setminus C = \overline{C} \union \{0\}
$.
From \eqref{eq: AxiomKPartOne}, 
$
G \intersect \overline{C} = \emptyset
$
since $0 \notin \overline{C}$.
Thus $\overline{C} \subset A \setminus \{0\}$ 
    satisfies $\AxiomK(H, G)$.
\end{proof}

A homomorphism $f : A \to A'$ of abelian groups is an epimorphism if $f$ is surjective. We use the notation $f : A \epimorphism A'$ to denote that $f$ is an epimorphism.

\begin{lemma} \label{lem: inductionStepForGraphIsNotConnected}
Let $C \subset A \setminus \{0\}$.
%, where $A$ is an abelian group. 
Let $B \le A$ be a subgroup that contains $C$.
Suppose that $C \subset B \setminus \{0\}$ satisfies $\AxiomK(G, H)$ for $G, H \le B$. 
%for some direct-sum decomposition $B = G \directSum H$ of $B$ into subgroups.
If the canonical epimorphism $A \epimorphism A/B$
    admits a section $s : A/B \to A$,
then $C \subset A \setminus \{0\}$ satisfies $\AxiomK(G, H \directSum s(A/B))$.  
\end{lemma}

\begin{proof}
Suppose that $s : A/B \to A$ is a section of $A \epimorphism A/B$
    and write $s(A/B) =: F$. 
Then $A = B \directSum F = (G \directSum H) \directSum F \isom G \directSum (H \directSum F)$. 
Since $C \subset B \setminus \{0\}$ satisfies $\AxiomK(G, H)$, so 
    \eqref{eq: AxiomKPartOne} holds in particular.
Next, since $C$ is contained in the summand $B$ of $A = B \directSum F$, so $H \subset B$ implies that 
$
(H \directSum F) \intersect C = H \intersect C = \emptyset,
$
using \eqref{eq: AxiomKPartTwo}. 
Thus $C \subset A \setminus \{0\}$ satisfies $\AxiomK(G, H \directSum F)$.
\end{proof}

Since exact sequences of vector spaces over any field 
    always split,
we have the following corollary.

\begin{corollary} \label{cor: inductionStepForGraphIsNotConnected}
Let $C \subset \K^d \setminus \{0\}$, where $\K$ is a field and $d \ge 0$ is an integer.
Suppose that $C \subset \langle C \rangle \setminus \{0\}$ satisfies $\AxiomK(V, W)$ for some vector subspaces $V, W \subset \langle C \rangle$. 
Then
$C \subset A \setminus \{0\}$ satisfies $\AxiomK(V, W \directSum s(\K^d/\langle C \rangle))$ for any section $s : \K^d/\langle C \rangle \to \K^d$
    of % the canonical $\K$-linear projection 
    $\K^d \epimorphism \K^d/\langle C \rangle$.  
\end{corollary}

\begin{proof}
Let $B := \langle C \rangle$ 
%to be the $\K$-linear span of $C$
    and apply Lemma \ref{lem: inductionStepForGraphIsNotConnected}.
\end{proof}

\section{Proof of Theorem \ref{thm: completeCores}} \label{sec: proofOfTheorem}

\subsection*{When $p = 2$:}

%This section is devoted to the proof of Theorem \ref{thm: completeCores} for $p = 2$, i.e.\ for cubelike graphs.

Let $d \ge 0$ be an integer throughout this section.

\begin{lemma} \label{lem: spanningSetsOfLowCardinalitySatisfyAxiomK}
Let $C \subset {\F_2}^d \setminus \{0\}$.
%, where $d \ge 0$ is an integer. 
If $\langle C \rangle_{\F_2} = {\F_2}^d$ 
 and $|C| < 5$,
then $C$ satisfies $\AxiomK(V, W)$
    for some vector subspaces $V, W \subset {\F_2}^d$.
\end{lemma}

\begin{proof}
Our proof is by exhaustion. See Table \ref{tab: cubelikeSmallCardinalitySpanning}.
Its third column gives an instance of $(V, W)$
    such that $C \subset {\F_2}^d \setminus \{0\}$ satisfies $\AxiomK(V, W)$.
\begin{center}
\begin{table}[h]
%\small
\begin{tabular}{| r |l | l | l|}
\hline
$d$ & $C \subset {\F_2}^d \setminus \{0\}$ & $(V, W)$ & $\dim_{\F_2}(V)$ \\ 
\hhline{|=|=|=|=|}
0 & $\emptyset $ & $(0, 0)$ & $0$ \\  
\hline 
1 & $\{i\} $ & $(\langle i \rangle, 0)$ & $0$ \\  
\hline
2 & $\{i, j\} $ & $(\langle i \rangle, \langle i + j\rangle)$ & $1$ \\  
& $\{i, j, i + j\} $ & $(\langle i, j \rangle, 0)$ & $2$ \\  
\hline
3 & $\{i, j, k\} $ & $(\langle i \rangle, \{\alpha i + \beta j + \gamma k: \, \alpha + \beta + \gamma = 0\})$ & $1$ \\  
& $\{i, j, k, i + j\} $ & $(\langle i, j \rangle, \langle i + j + k \rangle)$ & $2$ \\  
& $\{i, j, k, i + j + k\} $ & $(\langle i \rangle, \{\alpha i + \beta j + \gamma k: \, \alpha + \beta + \gamma = 0\})$ & $1$ \\ 
\hline
4 & $\{i, j, k, \ell\} $ & $(\langle i \rangle, \{\alpha i + \beta j + \gamma k + \delta \ell: \, \alpha + \beta + \gamma + \delta = 0\})$ & $1$ \\  
\hline
\end{tabular}
   \caption{Satisfaction of $\AxiomK$ by connection sets of small cardinality that span ${\F_2}^d$}
    \label{tab: cubelikeSmallCardinalitySpanning}   
\end{table}
\end{center} 
\vspace*{-\baselineskip}
%Here, we write $i := e_1$, $j := e_2$, $k := e_3$, and $\ell := e_4$,
%    where $\{e_1, \dots, e_d\}$ is a basis of ${\F_2}^d$. 
\end{proof}

\begin{proposition} \label{prop: setsOfExtremeCardinalitySatisfyAxiomK}
Let $C \subset {\F_2}^d \setminus \{0\}$.
%, where $d \ge 0$ is an integer. 
If either $|C| < 5$ or $|C| \ge 2^d - 5$,
then $C$ satisfies $\AxiomK(V, W)$
    for some vector subspaces $V, W \subset {\F_2}^d$.
\end{proposition}

\begin{proof}
Applying Corollary \ref{cor: inductionStepForGraphIsNotConnected}
    to Lemma \ref{lem: spanningSetsOfLowCardinalitySatisfyAxiomK}
        shows that every $C \subset {\F_2}^d \setminus \{0\}$
            with $|C| < 5$ satisfies $\AxiomK$.
Then, an application of Lemma \ref{lem: complementAlsoSatisfiesAxiomD}
    shows that every connection set in ${\F_2}^d \setminus \{0\}$
        with cardinality at least $2^d - 5$ also satisfies $\AxiomK$,
            since its complement in ${\F_2}^d \setminus \{0\}$ has cardinality at most $(2^d - 1) - (2^d - 5) = 4 < 5$.
        
Explicitly, the second column of Table \ref{tab: p=2ExtremeCardinality}
    gives an instance of $(V, W)$
    such that $C \subset {\F_2}^d \setminus \{0\}$ 
        with either $|C| < 5$ or $|C| \ge 2^d - 5$ 
    satisfies $\AxiomK(V, W)$.
The horizontal line in the middle of the table
    divides the connection sets $C$ into two,
        with $|C| < 5$ above the line
            and $|C| \ge 2^d - 5$ below the line.
In each row is an instance of $C$,
    where we write $s$ for any choice of section of the respective canonical linear projection.
For example, in the $4$th row where $C = \{i, j, i + j\}$,
    whose span is $\langle C \rangle = \langle i, j \rangle$,
        we write $s({\F_2}^d/\langle i, j \rangle)$
    to denote the image of a choice of section $s : {\F_2}^d/\langle i, j \rangle \to {\F_2}^d$
        of %the canonical $\F_2$-linear projection 
        ${\F_2}^d \epimorphism {\F_2}^d/\langle C \rangle =  {\F_2}^d/\langle i, j \rangle$,
            where $\{i, j\}$ is a linearly independent subset of ${\F_2}^d$.
             % Set some new page margins:
%\newgeometry{a4paper,left=1in,right=1in,top=1.2in,bottom=1.2in,nohead}
%\begin{landscape}
\begin{center}
\begin{table}[h]
\small
\begin{tabular}{| l | l | l|}
\hline
$C \subset {\F_2}^d \setminus \{0\}$ & $(V, W)$ & $\dim(V)$ \\ 
\hhline{|=|=|=|}
$\emptyset $ & $(0, {\F_2}^d)$ & $0$ \\  

$\{i\} $ & $(\langle i \rangle, s({\F_2}^d/\langle i \rangle))$ & $0$ \\  

$\{i, j\} $ & $(\langle i \rangle, \langle i + j \rangle \directSum s({\F_2}^d/(\langle i, j\rangle)))$ & $1$ \\  
$\{i, j, i + j\} $ & $(\langle i , j \rangle, s({\F_2}^d/\langle i , j \rangle))$ & $2$ \\  

$\{i, j, k\} $ & $(\langle i \rangle, \{\alpha i + \beta j + \gamma k: \, \alpha + \beta + \gamma = 0\} \directSum s({\F_2}^d/\langle i, j, k \rangle))$ & $1$ \\  
$\{i, j, k, i + j\} $ & $(\langle i, j \rangle, \langle i + j + k \rangle \directSum s({\F_2}^d/\langle i, j, k \rangle))$ & $2$ \\  
$\{i, j, k, i + j + k\} $ & $(\langle i \rangle, \{\alpha i + \beta j + \gamma k: \, \alpha + \beta + \gamma = 0\} \directSum s({\F_2}^d/\langle i, j, k \rangle))$ & $1$ \\ 

$\{i, j, k, \ell\} $ & $(\langle i \rangle, \{\alpha i + \beta j + \gamma k + \delta \ell: \, \alpha + \beta + \gamma + \delta = 0\} \directSum s({\F_2}^d/\langle i, j, k, \ell \rangle))$ & $1$ \\  
\hline
${\F_2}^d \setminus \{0\} $ & $({\F_2}^d, 0)$ & $d$ \\  

$\overline{\{i\}} $ & $(s({\F_2}^d/\langle i \rangle), \langle i \rangle)$ & $d$ \\  

$\overline{\{i, j\}} $ & $(\langle i + j \rangle \directSum s({\F_2}^d/(\langle i, j\rangle)), \langle i \rangle)$ & $d - 1$ \\  
$\overline{\{i, j, i + j\}}$ & $(s({\F_2}^d/\langle i , j \rangle), \langle i , j \rangle)$ & $d - 2$ \\  

$\overline{\{i, j, k\}}$ & $(\{\alpha i + \beta j + \gamma k: \, \alpha + \beta + \gamma = 0\} \directSum s({\F_2}^d/\langle i, j, k \rangle), \langle i \rangle)$ & $d - 1$ \\  
$\overline{\{i, j, k, i + j\}}$ & $(\langle i + j + k \rangle \directSum s({\F_2}^d/\langle i, j, k \rangle), \langle i, j \rangle)$ & $d - 2$ \\  
$\overline{\{i, j, k, i + j + k\}}$ & $(\{\alpha i + \beta j + \gamma k: \, \alpha + \beta + \gamma = 0\} \directSum s({\F_2}^d/\langle i, j, k \rangle), \langle i \rangle)$ & $d - 1$ \\ 
$\overline{\{i, j, k, \ell\}}$ & $(\{\alpha i + \beta j + \gamma k + \delta \ell: \, \alpha + \beta + \gamma + \delta = 0\} \directSum s({\F_2}^d/\langle i, j, k, \ell \rangle), \langle i \rangle)$ & $d - 1$ \\  
\hline
\end{tabular}
   \caption{Satisfaction of $\AxiomK$ by connection sets of extreme cardinality in ${\F_2}^d$}
    \label{tab: p=2ExtremeCardinality}   
\end{table}
\end{center}
%\end{landscape}
%\restoregeometry
\vspace*{-\baselineskip}

\end{proof} 

\begin{remark}
The bounds in Proposition \ref{prop: setsOfExtremeCardinalitySatisfyAxiomK} on the cardinality of $C$ are sharp as sufficient conditions for $C$ to satisfy $\AxiomK$. Specifically, when $d = 4$, there is a $C_* \subset {\F_2}^4 \setminus \{0\}$ with $|C_*| = 5$ such that neither $C_*$ nor its complement $\overline{C_*}$ satisfy $\AxiomK$. Note that, then $|\overline{C_*}| = 2^4 - 5 - 1$.
Indeed, we may take
\begin{equation}
C_* = \{i, j, k, \ell, i + j + k + \ell\} \subset {\F_2}^4 \setminus \{0\}.
\end{equation}  
Since $C_* \union \{0\}$ does not contain any $2$-dimensional vector subspace, thus any potential $V \subset C_* \union \{0\}$ is at most $1$-dimensional. Thus $W$ consists of all vectors $\alpha i + \beta j + \gamma k + \delta \ell$ that satisfy an equation $x\alpha + y\beta + z\gamma + w\delta = 0 $ for some scalars $x,y,z,w \in \F_2$.
As $i \notin W$, in order for $C_* \intersect W = \emptyset$, 
    thus $x = 1$. Similarly, $j, k, \ell \notin W$ imply that $y = z = w = 1$.
Thus the defining equation of $W$ is necessarily just $\alpha + \beta+ \gamma + \delta = 0$.
But this equation is satisfied by $i + j + k + \ell \in C_*$, where $(\alpha, \beta, \gamma, \delta) = (1, 1, 1, 1)$, so $i + j + k + \ell \in C_* \intersect W$, contradicting $C_* \intersect W = \emptyset$.
Therefore $C_*$ does not satisfy $\AxiomK$.

By Lemma \ref{lem: complementAlsoSatisfiesAxiomD}, $\overline{C_*}$ also does not satisfy $\AxiomK$. \qed
\end{remark}

%We are ready to prove Theorem \ref{thm: completeCores} in the case when $p = 2$.

\begin{proof}[Proof of Theorem \ref{thm: completeCores} when $p = 2$]
Apply Proposition \ref{prop: applicationOfAxiom}
    to Proposition \ref{prop: setsOfExtremeCardinalitySatisfyAxiomK}. 
%Note, in particular, as claimed in \S\ref{sec: mainResult}, if a cubelike graph $X$ on $({\F_2}^d, +)$ has degree $\deg(X) \ge 2^d - 5$,
%    then from Table \ref{tab: cubelikeExtremeCardinality},
%its core $X^\bullet = K_V$ has dimension at least $d - 2$.
\end{proof}

\subsection*{When $p = 3$:} %\label{sec: ternaryCayleyGraphs}

As the proof is entirely similar to the case when $p = 2$, we simply present our final solution and omit the arguments.

Before we present the proof, observe that connection sets in ${\F_3}^d \setminus \{0\}$ can be specified by their image in projective space ${\F_3} \P^{d - 1}$ under the canonical projection map 
$\pi : {\F_3}^d \setminus \{0\} \to  {\F_3} \P^{d - 1}$ given by $\pi(v) = \langle v \rangle_{\F_3} \setminus \{0\} = \{v, -v\} =: [v]$ for $v \in {\F_3}^d \setminus \{0\}$.
%When $d = 0$, regard ${\F_3}^0 \setminus \{0\} = \F_3 \P^{-1} = \emptyset$, so that $\pi = \mathrm{id}_\emptyset$ is just the identity map.
%The multiplicative group of our base field is ${\F_3}^\times = \{1, -1\}$.
Indeed $C = \pi^{-1}(\pi(C))$ because $-C = C$,
therefore $\pi(C) \subset {\F_3} \P^{d - 1}$ fully determines $C \subset {\F_3}^d \setminus \{0\}$.
%Hence, for every connection set , its defining property $-C = C$ implies that $\pi^{-1}(\pi(v)) = \{v, -v\} \subset C$ for every $v \in C$, so that
%\begin{equation}
%C = \pi^{-1}(\pi(C)).
%\end{equation}
%Therefore, to specify $C$,
%    it suffices to specify its image $\pi(C)$.
We thank xxxxxxxx (user:252194) for making this observation on Math Stack Exchange. %that connection sets in ${\F_3}^d \setminus \{0\}$ correspond bijectively to subsets of ${\F_3} \P^{d - 1}$ over Math Stack Exchange. 
(Their answer has been deleted, so we are unable to give a reference.)

\begin{proposition} \label{prop: setsOfExtremeCardinalitySatisfyAxiomKTernary}
Every $C \subset {\F_3}^d \setminus \{0\}$
%, where $d \ge 0$ is an integer. 
with $|C| < 12$ or $|C| \ge 3^d - 12$ satisfies $\AxiomK$. %direct-sum decomposition ${\F_3}^d = V \directSum W$
     %   of ${\F_3}^d$ into $\F_3$-vector subspaces.
\end{proposition}

\begin{proof}
Our proof is by exhaustion. We split our solution into two tables.
For $|C| < 12$ (resp.\ $|C| \ge 3^d - 12$), the fourth column of Table \ref{tab: p=3LowCardinality} (resp.\ the second column of Table \ref{tab: p=3HighCardinality}) gives an instance of $(V, W)$
    such that $C \subset {\F_3}^d \setminus \{0\}$ satisfies $\AxiomK(V, W)$.

 % Set some new page margins:
\newgeometry{a4paper,left=1.2in,right=1.2in,top=1.3in,bottom=1.3in,nohead}
\begin{landscape}   
\begin{center}
\begin{table}[h]
\small
\begin{tabular}{| r |c | l | l | l|}
\hline
$\dim(\langle C\rangle)$ & $|C|$ & $\pi(C) \subset \F_3 \P^{d- 1}$ & $(V, W)$ & $\dim(V)$ \\ 
\hhline{|=|=|=|=|=|}
0 & 0 & $\emptyset $ & $(0, {\F_3}^d)$ & $0$ \\ \hline 
1 & $2$ & $\{[i]\} $ & $(\langle i \rangle, s({\F_3}^d/\langle i \rangle))$ & $1$ \\  
\hline
2 & $4$ & $\{[i], [j]\} $ & $(\langle i \rangle, \langle i + j \rangle \directSum s({\F_3}^d/\langle i, j\rangle))$ & $1$ \\  
\hhline{|~|-|-|-|-|}
  & $6$ & $\{[i], [j], [i + j]\} $ & $(\langle i \rangle, \langle i - j \rangle \directSum s({\F_3}^d/\langle i, j \rangle))$ & $1$ \\ 
  \hhline{|~|-|-|-|-|}
  & $8$ & $\{[i], [j], [i + j], [i - j]\} $ & $(\langle i,  j \rangle,  s({\F_3}^d/\langle i, j \rangle))$ & $2$ \\ 
\hline
3 & $6$ & $\{[i], [j], [k]\} $ & $(\langle i \rangle, \{\alpha i + \beta j + \gamma k: \, \alpha + \beta + \gamma = 0\} \directSum s({\F_3}^d/\langle i, j, k \rangle))$ & $1$ \\  
\hhline{|~|-|-|-|-|}
  & $8$ & $\{[i], [j], [k], [i + j]\} $ & $(\langle k \rangle, \{\alpha i + \beta j + \gamma k: \, \alpha + \beta + \gamma = 0\} \directSum s({\F_3}^d/\langle i, j, k\rangle))$ & $1$ \\  
& & $\{[i], [j], [k], [i + j + k]\} $ & $(\langle i + j + k \rangle, \{\alpha i + \beta j + \gamma k: \, \alpha + \beta - \gamma = 0\} \directSum s({\F_3}^d/\langle i, j, k\rangle))$ & $1$ \\ 
\hhline{|~|-|-|-|-|}
  & $10$ & $\{[i], [j], [k], [j + k], [j - k]\} $ & $(\langle j, k \rangle, \langle i + j + k \rangle \directSum s({\F_3}^d/\langle i, j, k\rangle))$ & $2$ \\  
& & $\{[i], [j], [k], [i + k], [j + k]\} $ & $(\langle i \rangle, \{\alpha i + \beta j + \gamma k: \, \alpha + \beta + \gamma = 0\} \directSum s({\F_3}^d/\langle i,  j , k \rangle))$ & $1$ \\ 
& & $\{[i], [j], [k], [i + j - k], [j + k]\} $ & $(\langle i \rangle, \{\alpha i + \beta j + \gamma k: \, \alpha + \beta + \gamma = 0\} \directSum s({\F_3}^d/\langle i, j , k \rangle))$ & $1$ \\
\hline
4 & $8$ & $\{[i], [j], [k], [\ell]\} $ & $(\langle i \rangle, \{\alpha i + \beta j + \gamma k + \delta \ell: \, \alpha + \beta + \gamma + \delta = 0\} \directSum s({\F_3}^d/\langle i, j , k , \ell \rangle))$ & $1$ \\ 
\hhline{|~|-|-|-|-|}
  & $10$ & $\{[i], [j], [k], [\ell], [i + j]\} $ & $(\langle i + j \rangle, \{\alpha i + \beta j + \gamma k + \delta \ell: \, \alpha + \beta + \gamma + \delta = 0\} \directSum s({\F_3}^d/\langle i, j , k , \ell \rangle))$ & $1$ \\
  &  & $\{[i], [j], [k], [\ell], [i + j + k]\} $ & $(\langle \ell \rangle, \{\alpha i + \beta j + \gamma k + \delta \ell: \, \alpha + \beta -\gamma + \delta = 0\} \directSum s({\F_3}^d/\langle i , j , k , \ell \rangle))$ & $1$\\
  &  & $\{[i], [j], [k], [\ell], [i + j + k + \ell]\} $ & $(\langle i + j + k + \ell \rangle, \{\alpha i + \beta j + \gamma k + \delta \ell: \, \alpha + \beta + \gamma + \delta = 0\} \directSum s({\F_3}^d/\langle i, j , k , \ell \rangle))$ & $1$\\
\hline
5 & $10$ & $\{[i], [j], [k], [\ell], [m]\} $ & $(\langle i \rangle, \{\alpha i + \beta j + \gamma k + \delta \ell + \varepsilon m: \, \alpha + \beta + \gamma + \delta + \varepsilon = 0\} \directSum s({\F_3}^d/\langle i , j , k , \ell , m \rangle))$ & $1$ \\
\hhline{|=|=|=|=|=|}
\end{tabular}
   \caption{Satisfaction of $\AxiomK$ by connection sets of low cardinality in ${\F_3}^d$}
    \label{tab: p=3LowCardinality}   
\end{table}
\end{center}

\begin{center}
\begin{table}[h]
\small
\begin{tabular}{|l | l | l|}
\hline
$\pi(C) \subset \F_3 \P^{d- 1}$ & $(V, W)$ & $\dim(V)$ \\ 
\hhline{|=|=|=|}
$\F_3 \P^{d - 1}$ & $({\F_3}^d, 0)$ & $d$ \\ 
$\overline{\{[i]\}} $ & $(s({\F_3}^d/\langle i\rangle), \langle i \rangle)$ & $d - 1$ \\  
$\overline{\{[i], [j]\}} $ & $(\langle i + j \rangle \directSum s({\F_3}^d/\langle i, j \rangle), \langle i \rangle)$ & $d - 1$ \\  
$\overline{\{[i], [j], [i + j]\}} $ & $(\langle i - j \rangle \directSum s({\F_3}^d/(\langle i, j \rangle), \langle i \rangle)$ & $d - 1$ \\ 
$\overline{\{[i], [j], [i + j], [i - j]\}} $ & $(s({\F_3}^d/\langle i , j \rangle), \langle i,  j \rangle)$ & $d - 2$ \\ 
$\overline{\{[i], [j], [k]\}} $ & $(\{\alpha i + \beta j + \gamma k: \, \alpha + \beta + \gamma = 0\} \directSum s({\F_3}^d/\langle i , j, k \rangle), \langle i \rangle)$ & $d - 1$ \\  
$\overline{\{[i], [j], [k], [i + j]\}} $ & $(\{\alpha i + \beta j + \gamma k: \, \alpha + \beta + \gamma = 0\} \directSum s({\F_3}^d/\langle i, j, k \rangle), \langle k \rangle)$ & $d - 1$ \\  
$\overline{\{[i], [j], [k], [i + j + k]\}} $ & $(\{\alpha i + \beta j + \gamma k: \, \alpha + \beta -\gamma = 0\} \directSum s({\F_3}^d/\langle i,j, k\rangle), \langle i + j + k \rangle)$ & $d - 1$ \\ 
$\overline{\{[i], [j], [k], [j + k], [j - k]\}} $ & $(\langle i + j + k \rangle \directSum s({\F_3}^d/\langle i , j , k \rangle), \langle j , k \rangle)$ & $d - 2$ \\  
$\overline{\{[i], [j], [k], [i + k], [j + k]\}} $ & $(\{\alpha i + \beta j + \gamma k: \, \alpha + \beta + \gamma = 0\} \directSum s({\F_3}^d/\langle i , j , k \rangle), \langle i \rangle)$ & $d - 1$ \\ 
$\overline{\{[i], [j], [k], [i + j - k], [j + k]\}} $ & $(\{\alpha i + \beta j + \gamma k: \, \alpha + \beta + \gamma = 0\} \directSum s({\F_3}^d/\langle i , j , k \rangle), \langle i \rangle)$ & $d - 1$ \\
$\overline{\{[i], [j], [k], [\ell]\}} $ & $(\{\alpha i + \beta j + \gamma k + \delta \ell: \, \alpha + \beta + \gamma + \delta = 0\} \directSum s({\F_3}^d/\langle i , j , k , \ell \rangle), \langle i \rangle)$ & $d - 1$ \\ 
$\overline{\{[i], [j], [k], [\ell], [i + j]\}} $ & $(\{\alpha i + \beta j + \gamma k + \delta \ell: \, \alpha + \beta + \gamma + \delta = 0\} \directSum s({\F_3}^d/\langle i , j , k , \ell \rangle), \langle i + j \rangle)$ & $d - 1$ \\
$\overline{\{[i], [j], [k], [\ell], [i + j + k]\}} $ & $(\{\alpha i + \beta j + \gamma k + \delta \ell: \, \alpha + \beta -\gamma + \delta = 0\} \directSum s({\F_3}^d/\langle i , j , k , \ell \rangle), \langle \ell \rangle)$ & $d - 1$\\
$\overline{\{[i], [j], [k], [\ell], [i + j + k + \ell]\}} $ & $(\{\alpha i + \beta j + \gamma k + \delta \ell: \, \alpha + \beta + \gamma + \delta = 0\} \directSum s({\F_3}^d/\langle i , j , k , \ell \rangle), \langle i + j + k + \ell \rangle)$ & $d - 1$\\
 $\overline{\{[i], [j], [k], [\ell], [m]\}} $ & $(\{\alpha i + \beta j + \gamma k + \delta \ell + \varepsilon m: \, \alpha + \beta + \gamma + \delta + \varepsilon = 0\} \directSum s({\F_3}^d/\langle i , j , k , \ell , m \rangle), \langle i \rangle)$ & $d - 1$ \\
\hline
\end{tabular}
   \caption{Satisfaction of $\AxiomK$ by connection sets of high cardinality in ${\F_3}^d$}
    \label{tab: p=3HighCardinality}   
\end{table}
\end{center}
\end{landscape}
\restoregeometry % Restore the global document page margins
%Here $i := e_1$, $j := e_2$, $k := e_3$, $\ell := e_4$, and $m := e_5$,
%    where $\{e_1, \dots, e_d\}$ is a basis of ${\F_3}^d$.  
\end{proof}

%We are ready to prove our main theorem on tenary Cayley graphs, Theorem \ref{thm: TernaryCCForHighDegree}.

\begin{proof}[Proof of Theorem \ref{thm: completeCores} when $p = 3$]
Apply Proposition \ref{prop: applicationOfAxiom}
    to Proposition \ref{prop: setsOfExtremeCardinalitySatisfyAxiomKTernary}. 
\end{proof}

\subsection*{When $p \ge 5$:} %\ref{sec: largePrimesCase}

Let $p \ge 5$ be prime. 
If $C \subset {\F_p}^d \setminus \{0\}$ satisfies $|C| < \kappa(p) = 2$, then $|C| = 0$ so that $C = \emptyset$.
On the other hand, if $|C| \ge p^d - \kappa(p) = p^d - 2$,
then $|C| = p^d - 1$ so that $C = {\F_p}^d \setminus \{0\}$.
In either case, $C$ satisfies  $\AxiomK$:
%    for some direct-sum decomposition ${\F_p}^d = V \directSum W$
%        of ${\F_p}^d$ into $\F_p$-vector subspaces:

\begin{center}
\begin{table}[h]
%\footnotesize
\begin{tabular}{| l | l | l|}
\hline
$C \subset {\F_p}^d \setminus \{0\}$ & $(V, W)$ & $\dim_{\F_p}(V)$ \\ 
\hhline{|=|=|=|}
$\emptyset $ & $(0, {\F_p}^d)$ & $0$ \\  
\hline
${\F_p}^d \setminus \{0\} $ & $({\F_p}^d, 0)$ & $d$ \\  
\hline
\end{tabular}
   \caption{Satisfaction of $\AxiomK$ by connection sets of extreme cardinality in ${\F_p}^d$ for $p \ge 5$}
    \label{tab: p>=5ExtremeCardinality}   
\end{table}
\end{center}
\vspace*{-\baselineskip}
Thus Theorem \ref{thm: completeCores} for $p \ge 5$
    follows from Proposition \ref{prop: applicationOfAxiom}. \qed

\section{Details of Example \ref{eg: counterexamples}} \label{sec: counterexamples}

We provide explicit constructions below.
The second column is the rank $d_\ast$ of the abelian group ${\Z_p}^{d_*}$ that the counterexample $X_* = \Cay({\Z_p}^{d_\ast}, C_*)$ is a Cayley graph on.
%We write $i := e_1$, $j := e_2$, $k := e_3$, and $\ell := e_4$,
%where $\{e_1, \dots, e_{d_\ast}\}$ is a $\Z_p$-linear basis of ${\Z_p}^{d_\ast}$.
\begin{center}
\begin{table}[h]
\small
\begin{tabular}{| r ||l | l  l|}
\hline
$p$ & $d_\ast$ & $C_\ast \subset {\Z_p}^{d_\ast} \setminus \{0\}$ & \textbf{Name of } $X_\ast := \Cay({\Z_p}^{d_\ast}, C_\ast)$ \\ 
\hhline{|=#=|= =|}
& & &\\[-1.3em]
%\rule{0pt}{0.6cm}
$2$ & $4$ & $\{i, j, k, \ell, \ijkl\}$ & Folded $5$-cube/ Clebsch graph %$\overline{\frac{1}{2}Q_5}$ 
\\
% & & $C^* = \{i, j, k, \ell, i + j, i + k, i + \ell, j + k, j + \ell, k + \ell\}$ & Halved $5$-cube $\frac{1}{2}Q_5$\\  
\hline 
$3$ & $3$ & $\pi^{-1}(\{[i], [j], [k], [i + j], [i + k], [i + j + k]\})$ & (No standard name) \\  
% &  & $C^* = $ & No name \\
\hline
$\ge 5$ & $1$ & $\{i, -i\}$ & $p$-cycle $C_p$ \\  
% &  & $C^* = \F_p \setminus \{0, i, -i\}$ & $\overline{C_p}$ \\ 
\hline
\end{tabular}
   \caption{$p$-ary Cayley graph cores of degree $\kappa(p)$ whose complements are also cores}
    \label{tab: counterexamples}   
\end{table}
\end{center}
\vspace*{-\baselineskip}
Note that $\overline{X_\ast} = \Cay({\Z_p}^d, \overline{C_\ast})$. Also note that neither $X_\ast$ nor $\overline{X_\ast}$ is complete.

We proceed to give brief explanations or citations to the literature why $X_\ast$ and $\overline{X_\ast}$ are cores.

\paragraph{For $p = 2$:} Each endomorphism $f$ of the folded $5$-cube $X_\ast$ sends every $5$-cycle to another $5$-cycle since $X_\ast$ has no triangles. In particular, $f$ is injective on every $5$-cycle. Any two vertices of $X_\ast$ lie on a common 5-cycle, so $f$ is globally injective, and therefore invertible. Therefore $X_\ast$ is a core. This proof is due to Lavrov \cite{Lavrov2024}.
%, who answered to a question of the second author on Math Stack Exchange. 

Godsil and Royle \cite[Theorem 4.6]{GodsilRoyle2011} showed that the halved $n$-cube 
\begin{equation}\frac{1}{2}Q_n = \Cay({\Z_2}^{n - 1}, \{e_1, \dots, e_{n - 1}\} \union \{e_i + e_j :\, 1 \le i <j \le n - 1\})
\end{equation}
is a core if and only if $n$ is not a power of $2$. In particular, 
$\overline{X_\ast} \isom \frac{1}{2}Q_5$ is a core.

\paragraph{For $p = 3$:} Since $X_\ast$ has clique number $4$, hence ${X_\ast}^\bullet$ is not complete. Suppose on the contrary that $X_\ast$ is not a core. Then, by a result of Hahn and Tardif \cite{HahnTardif1997}, $|V({X_\ast}^\bullet)|$ divides $|V(X_\ast)| = 27 = 3^3$. Hence $|V({X_\ast}^\bullet)| = 9$ since $X^\bullet$ is not complete. By general theory, ${X_\ast}^\bullet$ is a $9$-vertex vertex-transitive graph with the same chromatic number as $X_\ast$. Since $X_\ast$ has chromatic number at least $7$, we obtain a contradiction as there are no non-complete $9$-vertex vertex-transitive graphs with chromatic number more than $5$. Therefore, $X_\ast$ is a core.

The proof that $\overline{X_\ast}$ is a core is exactly the same, using the fact that $\overline{X_\ast}$ has clique number $4$ and chromatic number at least $7$. These proofs for $p = 3$ are due to Royle \cite{Royle2024}.

\paragraph{For $p \ge 5$:} We prove more generally that any Cayley graph $X$ on $\Z_p$ of positive degree is a core. Indeed, by Hahn and Tardif's result again, either $|V(X^\bullet)| = 1$ or $|V(X^\bullet)| = p$ because $|V(X)| = p$ is prime. Assuming that $\deg(X) > 0$, the map that sends all the vertices of $X$ to a single vertex is not an endomorphism of $X$, thus $|V(X^\bullet)| \neq 1$. Therefore $|V(X^\bullet)| = p$, so that $X^\bullet = X$. In particular, $X_\ast = C_p$ and $\overline{X_\ast} = \overline{C_p}$ on $\Z_p$ are cores. Note that $\deg(\overline{X_\ast}) = p - 3 > 0$ since $p \ge 5$. 

%Finally, we remark that since $\deg(X_\ast) = \kappa(p) \neq p^d - 1$ for any integer $d \ge 0$, so $X_\ast$ is indeed not complete. Also $\kappa(p) > 0$, so $X_\ast$ is not edgeless, hence $\overline{X_\ast}$ is also not complete. 
\medskip
In summary, for each prime $p$, we have constructed non-complete $p$-ary Cayley graph cores $X_*$ and $\overline{X_\ast}$ of degree $\kappa(p)$ and $|V(\overline{X_\ast})| - \kappa(p) - 1$ respectively. Hence neither ${X_*}^\bullet = X_*$ nor ${\overline{X_\ast}}^\bullet = \overline{X_\ast}$ are complete. Therefore the degree bounds given in Theorem \ref{thm: completeCores} are optimal.

%\section{Conclusion} \label{sec: Conclusion}

%Man\v{c}inska, Pivotto, Roberson and Royle wrote the main substantial paper addressing Ne\v{s}et\v{r}il and \v{S}\'{a}mal's question that the core of a cubelike graph is cubelike \cite{MPRR2020}. They proved this conjecture in the affirmative for cubelike graphs with small cores, namely cores with at most $32$ vertices. In this paper, we gave a sufficient condition for the core of a cubelike graph to be complete and hence cubelike (see Proposition \ref{prop: applicationOfAxiom}). We showed that this sufficient condition is satisfied for cubelike graphs $X$ of extreme degree, i.e.\ whenever $\deg(X) < 5$ or $\deg(X) \ge |V(X)| - 5$ (see Proposition \ref{prop: setsOfExtremeCardinalitySatisfyAxiomK}). More specifically, when the degree of a cubelike graph $X$ is high, namely $\deg(X) \ge |V(X)| - 5$, then there is some $\F_2$-linear subspace of $V(X)$ of dimension at least $\dim_{\F_2}(V(X)) - 2$ that induces a core of $X$, see Table \ref{tab: cubelikeExtremeCardinality}. Letting $V(X) \to \infty$, cubelike graphs of high degree have arbitrary large cores, thereby answering the conjecture in the affirmative for a class of cubelike graphs distinct from that addressed by Man\v{c}inska, Pivotto, Roberson and Royle's result. 

\paragraph{Conclusion.} Although the counterexample $X_\ast$ (resp.\  $\overline{X_\ast}$) shows that not every $p$-ary Cayley graph of degree $\kappa(p)$ (resp.\ of degree $|V(\overline{X_\ast})| - \kappa(p) - 1$) has a complete core,
there remains this question:
\begin{question}
For fixed $d$, is there any integer $\kappa(p) + 1\le n \le p^d - \kappa(p) - 2$ such that every Cayley graph on ${\Z_p}^d$ of degree $n$ has a complete core? 
\end{question}
When $p = 2$, this question appears open for $d \ge 7$.
%We formulate the following question:
%\begin{description}
%\item[Question C:] Fix an integer $n \ge 0$. What is the largest integer $d_0 = d_0(n) \ge 0$, depending on $n$, such that every cubelike graph with $2^d$ vertices of degree $n$ has a complete core for all $d \le d_0$? 
%\end{description}

%For example, Theorem \ref{thm: completeCores} gives $d_0(0) = d_0(1) = d_0(2) = d_0(3) = d_0(4) = \infty$. Also Man\v{c}inska, Pivotto, Roberson and Royle's result \cite[Theorem 8.1]{MPRR2020} together with the counterexample $\overline{\frac{1}{2}Q_5}$ in Example \ref{eg: counterexamples} implies that $d_0(5) = 3$.

%Note that, if $d_0(n) < \infty$, then there is, by definition, some cubelike graph $X_n$ with $2^{d_0(n) + 1}$ vertices of degree $n$ whose core is not complete. Then, in fact, for every $d \ge d_0(n) + 1$, the disjoint union $\coprod_{i = 1}^{2^{d - d_0 - 1}} X_n$ of $2^{d - d_0 - 1}$ copies of $X_n$ is a cubelike graph with $2^d$ vertices of degree $n$ whose core $(\coprod_{i = 1}^{2^{d - d_0 - 1}} X_n)^\bullet = {X_n}^\bullet$ is not complete.
%Thus, the computation of $d_0(n)$ as a function of $n$ in Question C would give a complete solution of when every cubelike graph with $2^d$ vertices of degree $n$ has a complete core, namely exactly when $d \le d_0(n)$.

\paragraph{Acknowledgements.}
We would like to thank Professor Gordon Royle for reading our paper and providing helpful comments, observations and data via private communication. 
The second author asked several questions on Math Stack Exchange and Math Overflow
    and would like to thank
ancient mathematician,
ahulpke,
schiepy,
Steve D,
Ewan Delanoy,
Jyrki Lahtonen, 
Misha Lavrov,
Alex Ravsky, 
Gordon Royle, 
Richard Stanley, and
xxxxxxxxx
for their comments and answers.

This research was supported by the Yale-NUS College seed grant (IG23-SG009).

%===================================================================%

\end{document}